\documentstyle[11pt]{article}
\input{amssym.def}
\input epsf

\newtheorem{thm}{Theorem}[section]
\newtheorem{lemma}[thm]{Lemma}

\newtheorem{rmk}[thm]{Remark}

\newcommand{\B}{{\rm B}}
\newcommand{\C}{{\Bbb C}}

\newcommand{\Z}{{\Bbb Z}}

\newcommand{\cA}{{\cal A}}

\newcommand{\cL}{{\cal L}}

\newcommand{\cP}{{\cal P}}

\newcommand{\cU}{{\cal U}}

\newcommand{\cD}{{\cal D}}

\newcommand{\fm}{{\frak m}}

\begin{document}
\bibliographystyle{plain}

\title{Versal Deformations of Formal Arcs\footnote{AMS-MSC: 14B20
% ; this a is a draft, we will be grateful for any 
% comments or bibliographical suggestions.
}}
\author{Mikhail Grinberg and David Kazhdan}
% \date{October 22, 1998}
\maketitle

\section{Introduction}

Let $X$ be a complex algebraic variety, and $X^\circ \subset X$
be the smooth part of $X$.  Consider the scheme $\cL(X)$ of
formal arcs in $X$.  The $\C$-points of $\cL(X)$ are just maps
$D = {\rm Spec} \, \C[[t]] \to X$ (see, for example, [DL] for a
definition of $\cL(X)$ as a scheme).  Let $\cL^\circ(X)$ be the
open subscheme of arcs whose image is not contained in
$X \setminus X^\circ$.  Fix an arc $\gamma : D \to X$ in
$\cL^\circ(X)$, and let $\cL(X)_\gamma$ be the formal
neighborhood of $\gamma$ in $\cL(X)$.  The purpose of this paper
is to prove the following.

\begin{thm}\label{main}
There exists a scheme $Y = Y(\gamma)$ of finite type over $\C$,
and a point $y \in Y(\C)$, such that:
$$\cL(X)_\gamma \cong Y_y \times D^\infty,$$
the formal neighborhood of $y$ in $Y$, times a product of
countably many copies of $D$.
\end{thm}

The finite-dimensional piece $Y_y$ may be called the parameter
space of a versal deformation of $\gamma$.  It gives a model
for the singularity of $\cL(X)$ at $\gamma$.  It is not hard
to check that the analytic germ $(Y, y)$ is determined
uniquely up to multiplying by a finite-dimensional vector
space.  We note, as a pitfall, that the embedding $Y_y \to
\cL^\circ(X)$ given by Theorem \ref{main} does not extend to
a map of the analytic germ $(Y, y) \to (\cL^\circ(X), \gamma)$.
Intuitively, this is because the versal deformation captures
the breaking-up of singularities, but the formal disc has only
one point.

The study of formal arcs in an algebraic variety was originated
by J. Nash in [Na].  It has been a focus of much recent
activity (see [Ba], [DL]).  However, this recent work has
mostly concentrated on the spaces of truncated arcs.

We expect the singularities of the spaces $Y(\gamma)$ to arise
in many problems involving maps of smooth curves into singular
varieties, or other objects which can be locally described
as such maps.  The original motivating example of this
is the work of Feigin, Finkelberg, Kuznetsov, and Mirkovi\'c
(see [K], [FK], [FFKM], and their references) on quasimap
spaces $Q_\alpha^D$.  Conjecturally, the singularities of the
$Q_\alpha^D$ are the same as those of $\cL^\circ (X)$, where
$X = \overline{G/U}$ is the affine closure of the quotient
of an algebraic group by the unipotent radical of a Borel.

Theorem \ref{main} was conjectured by V. Drinfeld [Dr].
Drinfeld also proved a weaker version of it, in which the
finite-dimensional piece $Y_y$ is not identified as the
formal neighborhood of a point of a scheme of finite type.

The paper is organized as follows.  In Section 2 we prove a
lemma about plane curves depending on parameters on which the
proof of Theorem \ref{main} is based.  Theorem \ref{main} is
proved in Section 3, and Section 4 contains a few simple
examples.

We would like to thank V. Drinfeld for generously sharing his
ideas.  We are also grateful to A. Braverman and D. Gaitsgory
for many useful conversations, and to R. Vakil for his help
with Lemma \ref{standard}.

\vspace{.1in}

\noindent
{\bf A notational convention.}  Throughout this paper,
a {\it test-ring} $A$ is a finite-dimensional, local
commutative $\C$-algebra with $1$.  If $S$ is a scheme
over $\C$, and $p \in S(\C)$ is a $\C$-point, we think of the
formal neighborhood $S_p$ in terms of its functor of points:
$A \mapsto S_p(A)$, from test-rings to sets.

\section{Plane Curves Depending on Parameters}

In this section, we prove a lemma about plane curves
depending on parameters which will be used in the
proof of Theorem \ref{main}.  Let $f(x,y)$ be a
polynomial vanishing at the origin.  Denote the
curve $\{ f (x,y) = 0 \}$ by $C$.  Let $\cU$ be a
small neighborhood of $0$ in $\C^2$.  Assume that
$f$ has no critical points in $\cU \setminus \{ 0 \}$.
In addition, assume that $C \cap \cU = C_1 \cup C_2$, where
$C_1$ is a smooth analytic branch of $C$, which is tangent
to the $x$-axis at $0$, and $C_1 \cap C_2 = \{ 0 \}$.  In
other words, we have $f|_\cU = f_1 \cdot f_2$, where
$f_1, f_2 : \cU \to \C$ are analytic functions, the
differential $d_0 f_1 = y$, and $f_1 (x,y) = f_2(x,y) = 0$
implies $x = y = 0$.

Let $m$ be the multiplicity of $0$ as a point of
intersection of $C_1$ and $C_2$.  In other words, $m$
is the order of vanishing of the restriction $f_2 |_{C_1}$
at $0$.  Pick a large integer $M$, and let $\C[x,y]^{\leq M}
\subset \C[x,y]$ be the affine space of all polynomials
of degree $\leq M$.  Denote by $P$ a small analytic
neighborhood of $f$ in $\C[x,y]^{\leq M}$.  Let
$R \subset P$ be the locus of $g \in P$, such that the
curve $\{ g(x,y) = 0 \}$ has an analytic branch near the
origin which is $C^\infty$-close to $C_1$.  We may
choose a neighborhood $\cD$ of zero in $\C$, such that
for any $g \in R$, there is a unique analytic function
$h_g : \cD \to \C$ satisfying $g (t, h_g (t)) = 0$
for all $t \in \cD$.  Let $b_k (g)$ be the $k$-th coefficient
of the Taylor series for $h_g$, so that $h_g (t) =
\sum_{k = 0}^\infty b_k (g) \cdot t^k$.  Finally, let
$R_f \subset P_f$ be the formal neighborhoods of $f$ in $R$
and $P$.

\begin{lemma}\label{plane_curves}
(i)   $R$ is a smooth analytic subvariety of $P$, of
codimension $m$.

(ii)  $R$ is a local analytic branch of some algebraic
subvariety of $\C[x,y]^{\leq M}$.  

(iii) The functions $b_k : R \to \C$ are complex analytic.

(iv)  Let $A$ be a test-ring with maximal ideal $\fm$.
Then $R_f(A)$ is the set of all $\tilde f \subset
A [x,y]^{\leq M}$ such that $\tilde f \equiv f$
{\em (mod $\fm$)}, and there exists an $\tilde h \in A [[t]]$
satisfying $\tilde h \equiv h_f$ {\em (mod $\fm$)} and
$\tilde f (x, \tilde h (x)) \equiv 0$.
\end{lemma}

Our proof of Lemma \ref{plane_curves} will use the following
basic result of Tougeron (see [To]).

\begin{lemma}\label{technical} {\em [AVG, Part I, Ch 6.3]}
Let $f: (\C^r, 0) \to (\C, 0)$ be a germ of a complex
analytic function with an isolated singularity at the origin.
Then there exists an $n \in \Z_+$, such that any germ
$f': (\C^r, 0) \to (\C, 0)$, with the same $n$-jet as $f$, is
analytically right-equivalent to $f$.
\hfill $\Box$
\end{lemma}

Both [To] and [AVG] state Lemma \ref{technical} for real
$C^\infty$ functions.  However, the proof in [AVG] goes through
word-for-word in the complex analytic setting.

\vspace{.1in}

\noindent
{\bf Proof of Lemma \ref{plane_curves}:}
The set $P$ gives a right-versal deformation of the
singularity of $f$ at the origin (see [AVG] for a discussion
of versal deformations of singularities of functions).
Therefore, to prove (i), (iii), and (iv), we can replace $f$
by any polynomial $f'$ which satisfies the hypothesis of the
lemma, and is analytically right-equivalent to $f$ near zero.
Using Lemma \ref{technical}, it is not hard to check
that such an $f'$ can be chosen to be divisible by
$y$.  In other words, in the proof of parts (i), (iii), and
(iv), we can assume that $C_1$ is the $x$-axis.  We proceed
with this assumption.

Part (iii) is a combination of the Cauchy integral formula
and the implicit function theorem.  To be precise, fix a
small $\epsilon > 0$, and let $S_\epsilon$ be the
$\epsilon$-circle in the $x$-axis.  Then the $\{b_k (g)\}$ can
be computed as the positive Fourier coefficients of the
restriction $h_g |_{S_\epsilon}$.  But this restriction
depends analytically on $g$, by the implicit function theorem.

Note that this argument shows that we have a uniform bound on
the size of the $b_k = b_k (g)$ in terms of the size of $g - f$.
More precisely, let
$$g(x,y) = f(x,y) + \sum_{i,j} c_{i,j} \, x^i y^j$$
($0 \leq i,j \leq M$).  Then there exists a $\sigma > 0$ such
that:
\begin{equation}\label{bsmall}
|b_k| < \sigma \cdot \epsilon^{-k} \cdot \max_{i,j} |c_{i,j}|,
\end{equation}
for any $g \in R$, and any $k \in \Z_+$.

Turning now to part (i), let ${\displaystyle f(x,y) =
\sum_{i,j} a_{i,j} \, x^i y^j}$ ($a_{i,j} \in \C$). 
By assumption, $a_{i,0} = 0$ for all $i$, and $a_{i,1} = 0$ for
$i < m$.  Our proof of part (i) is based on analyzing the
identity $g(x, h_g(x)) = 0$.  Specifically, setting each
coefficient of the power series $g(x, h_g(x)) \in \C[[x]]$ to
zero, we obtain an infinite system of equations in the
$\{c_{i,j}\}$ and the $\{b_k\}$.  Let us inspect the first
few equations of this system.  We write $O(b^2)$ for terms that
contain a product of at least two of the $b_k$.

\vspace{.1in}

$$c_{0,0} + c_{0,1} \, b_0 + O(b^2) = 0$$
$$c_{1,0} + c_{1,1} \, b_0 + c_{0,1} \, b_1 + O(b^2) = 0$$
$$c_{2,0} + c_{2,1} \, b_0 + c_{1,1} \, b_1 +
c_{0,1} \, b_2 + O(b^2) = 0$$
\begin{equation}\label{system}
\dots
\end{equation}
$$c_{m,0} + a_{m,1} \, b_0 + c_{m,1} \, b_0 +
c_{m-1,1} \, b_1 + \dots + c_{0,1} \, b_m + O(b^2) = 0$$
$$c_{m+1,0} + a_{m+1,1} \, b_0 + c_{m+1,1} \, b_0 +
c_{m,1} \, b_1 + \dots + c_{0,1} \, b_{m+1} + O(b^2) = 0$$
$$\dots$$

\vspace{.1in}

We are interested in finding small solutions $\{c_{i,j}, b_k\}$
of the system (\ref{system}).  For this, let us first fix a
test-ring $A$ with maximal ideal $\fm$, and look for solutions
$\{c_{i,j}, b_k\}$ in $\fm$ (while the $a_{i,j}$ are still in $\C$).
Let $I = \{ 0, \dots, M\}^2 \setminus \{ 0, \dots, m-1\} \times
\{ 0 \}$.  Then it is obvious from inspecting the system
(\ref{system}) that it has a unique solution for any $\{c_{i,j} \in
\fm\}_{(i,j) \in I}$.  Indeed, we can first construct the solution
modulo $\fm^2$ by solving the first equation for $c_{0,0}$, the
second equation for $c_{1,0}$, .\,.\,.\,, the $m$-th equation for
$c_{m-1,0}$, the $(m+1)$-st equation for $b_0$, the $(m+2)$-nd
equation for $b_1$, and so on.  Note that at this stage we have
$c_{0,0}, \, c_{1,0}, \, \dots, \, c_{m-1, 0} \in \fm^2$.  Once the
solution is constructed modulo $\fm^2$, we can go back to the first
equation and solve for $c_{0,0}$ modulo $\fm^3$, and so on.

This observation means that there are power series
$s_0, \dots, s_{m-1} \in \C[[c_{i,j}]]_{(i,j) \in I}$, such
that for any solution $\{c_{i,j}, b_k\}$ of (\ref{system})
in the maximal ideal of any test-ring, we have:
\begin{equation}\label{R}
c_{l,0} = s_l (c_{i,j})_{(i,j) \in I}, \;\; \mbox{for} \;\;
l = 0, \dots, m-1.
\end{equation}
It is a straightforward exercise in combinatorics to show that the
coefficients of the power series $s_l$ can not grow
super-exponentially, and that the $s_l$ all converge in some
neighborhood of zero in $\C^I$.  It follows that the locus
$R \subset P$ is the submanifold defined by the equations (\ref{R}),
where the $c_{i,j}$ are now small complex numbers.  The passage from
the formal solution of (\ref{system}) to the solution in small
$\{c_{i,j}, \, b_k\}$ presents no difficulties because of the
estimate (\ref{bsmall}).  This completes the proof of (i).

For part (iv), note that
$$P_f (A) = \{ \tilde f \subset A [x,y]^{\leq M} \; | \;
\tilde f \equiv f \; \mbox{(mod $\fm$)} \}.$$
Any $\tilde f \in P_f (A)$ can be written as
$\displaystyle \tilde f (x,y) = f (x,y) +
\sum_{i,j} c_{i,j} \, x^i y^j$, with $c_{i,j} \in \fm$
($0 \leq i,j \leq M$).  The point $\tilde f$ is in $R_f(A)$ if and
only if the $\{ c_{i,j} \}$ satisfy equations (\ref{R}).  By
construction, this means that there exist
$\{ b_k \in \fm \}_{k \in \Z_+}$, such that $\{c_{i,j}, b_k\}$ is a
solution of (\ref{system}).  Setting $\tilde h = \sum_{k = 0}^\infty
b_k \cdot t^k$ proves part (iv).

For part (ii), let $Z \subset \C[x,y]^{\leq M}$ be the set of all
polynomials $g(x,y)$ such that the curve $\{ g(x,y) = 0 \}$
has at least $m$ singular points (counted with multiplicities).
Then $Z$ is a closed subvariety of $\C[x,y]^{\leq M}$, containing
$R$.  A generic point $g \in R$ corresponds to a curve with $m$
simple double points $q_1, \dots, q_m$ near the origin.  Let
$\tilde g \in P$ be a polynomial very near $g$.  For each of the
$q_i$, the condition on $\tilde g$ that the curve $\{ \tilde g (x, y)
= 0 \}$ has a double point near $q_i$ is smooth of codimension $1$.
Furthermore, if $M$ is sufficiently large, these conditions for
different $q_i$ are independent.  It follows that the codimension
of $Z$ at $g$ is $m$.  Therefore, $R$ is a local analytic branch of
$Z$.
\hfill$\Box$

\begin{rmk}
{\em  The proof above shows that the tangent space
$T_f R \subset T_f P$ is the coordinate plane
$c_{0,0} = c_{1,0} = \dots = c_{m-1,0} = 0$.}
\end{rmk}

\section{Proof of Theorem \ref{main}}

We break the proof up into six steps.

\vspace{.1in}

\noindent
{\bf Step 1: Approximating by analytic arcs.}
Without loss of generality, we may assume that $X$ is
a closed subvariety of a finite-dimensional vector space
$U$, and that $\gamma (0)$ is the origin.  Below, we state
three preliminary lemmas whose proofs are routine.
For $n \in \Z_+$, let $D_n = {\rm{Spec}} \, \C[t] / t^{n+1}$,
and  let $\cA(X, \gamma, n)$ be the set of all arcs
$\alpha : D \to X$ which agree with $\gamma$ on $D_n$.

\begin{lemma}\label{prelim1}
For any $n \in \Z_+$, there exists an analytic
$\alpha \in \cA(X, \gamma, n)$.
\hfill$\Box$
\end{lemma}

\begin{lemma}\label{prelim2}
There is an $n(\gamma) \in Z_+$, such that the limit
$\lim_{t \to 0} T_{\alpha(t)} X$ is the same for all analytic
$\alpha \in \cA(X, \gamma, n(\gamma))$.  We call this limit
$V = V(\gamma)$; it is a linear subspace of $U$ of dimension
$d = \dim X$.
\hfill$\Box$
\end{lemma}

Choose a direct sum decomposition $U = V \oplus W$. 
Let $p: U \to V$ be the projection along $W$, and
$p_* : \cL(X) \to \cL(V)$ be the map induced by $p$.
We will use the abbreviation $p_*\gamma = p\gamma$.
Denote by $\cA(V, p\gamma, n)$ the set of all arcs
$\beta : D \to V$ which agree with $p\gamma$ on $D_n$.

\begin{lemma}\label{prelim3}
The number $n(\gamma)$ in Lemma \ref{prelim2} can be chosen
so that, for any choice of $W$, the map $p_*$
induces a bijection $\cA(X, \gamma, n(\gamma)) \cong
\cA(V, p\gamma, n(\gamma))$ which takes analytic arcs
to analytic arcs.
% \hfill$\Box$
\end{lemma}

\noindent
{\bf Proof:}  Fix a Hermitian metric on $U$.  Let $\B(x, r)$
($x \in X$, $r > 0$) be the $r$-ball around $x$.  It suffices
to choose $n(\gamma)$ so that for some analytic $\alpha \in
\cA(X, \gamma, n(\gamma))$, and every small $t \neq 0$,
the restriction of $p$ to $X \cap \B(\alpha(t), \,
|t|^{n(\gamma) - 1})$ is injective.
\hfill$\Box$

\vspace{.1in}

\noindent
{\bf Step 2: Truncating the arc.}
Since $V$ is a vector space, an arc in $V$ has well-defined
higher derivatives.  Therefore, for any $N \in \Z_+$, we can
write
\begin{equation}\label{decomp}
\cL(V) = \cL(V)^{\leq N} \times \cL(V)^{>N},
\end{equation}
where $\cL(V)^{\leq N} = \{ \beta : D \to V \; | \;
\beta^{(n)} (0) = 0, \; \mbox{for} \; n > N \}$, and
$\cL(V)^{> N} = \{ \beta : D \to V \; | \;
\beta^{(n)} (0) = 0, \; \mbox{for} \; n \leq N \}$.

\begin{lemma}\label{truncate}
For any sufficiently large $N$, there exists an analytic
arc $\alpha : D \to X$ such that:

(i)   $\alpha$ agrees with $\gamma$ on $D_N$;

(ii)  $p_* \alpha \in \cL(V)^{\leq N}$; and

(iii) $\cL(X)_\alpha$ is isomorphic to $\cL(X)_\gamma$.
\end{lemma}

\noindent
{\bf Proof:}  The idea is to use the fact that $\cL (X)$
has a very rich set of symmetries.  Specifically, any map
$t \mapsto \xi_t$ from the formal disc to the vector fields
on $X$ produces a vector field $\hat \xi$ on $\cL (X)$.
Such a vector field, in turn, produces a time-$1$ flow map
$\Xi_1 : \cL (X)^\xi \to \cL (X)$.  This $\Xi_1$ should be
understood as follows.  Let $X^\xi(\C)$ be the set of
$\C$-points of $X$ where the time-$1$ flow of $\xi_0$ is
defined (it is an open subset of $X(\C)$ in the analytic
topology).  Let now $A$ be a test-ring.  Denote by
$\cL(X)^\xi (A) \subset \cL(X)(A)$ the set of all maps
${\rm Spec} \, A \times D \to X$ which send the unique
$\C$-point of the domain to $X^\xi(\C)$.  Then $\Xi_1$ is
a natural transformation between the functor $A \mapsto
\cL(X)^\xi (A)$ and  the functor $A \mapsto \cL(X)(A)$.

Pick a set of coordinates $(u_1, \dots, u_d)$ on $V$, and
consider the basic vector fields
$(\frac{\partial}{\partial u_1},
\dots, \frac{\partial}{\partial u_d})$.
They lift to meromorphic vector fields
$(p^*\frac{\partial}{\partial u_1}, \dots,
p^*\frac{\partial}{\partial u_d})$ on $X^\circ$.  We may
choose regular functions $(g_1, \dots, g_d)$ on $X$ to
`cancel out the poles' of $(p^*\frac{\partial}{\partial u_1},
\dots, p^*\frac{\partial}{\partial u_d})$, so that each
$\eta_i = g_i \cdot p^*\frac{\partial}{\partial u_i}$
is a regular vector field on $X$.  Note that the image of
$\gamma$ is not contained in the closure of the critical
set of the restriction of $p$ to $X^\circ$.  This means that
the $g_i$ can be chosen not to vanish identically along
$\gamma$.  Let $n_i$ be the order of vanishing of $g_i$ along
$\gamma$, that is, the smallest $n$ such that the $n$-th
coefficient $(g_i \circ \gamma)_n \neq 0$.  Set $N_0 = \max
(n_1, \dots, n_d, n(\gamma))$, where $n(\gamma)$ is as in
Step 1.  The claims of the lemma will hold for any $N \geq
N_0$.  Fix such a number $N$.  We construct the arc $\alpha$
in a series of $d$ steps.  As a first approximation we take
$\alpha_0 = \gamma$.  The first step is accomplished by the
following claim.

\vspace{.1in}

\noindent
{\bf Claim:}  There is a unique $f \in t^{N + 1 - n_1} \cdot
\C[[t]]$ with the following property.  Define a family $\xi_t$
of vector fields on $X$ by $\xi_t = f(t) \cdot \eta_1$.  Let
$\hat \xi$ be the corresponding vector field on $\cL (X)$, and
$\Xi_1$ be the time-$1$ flow of $\hat \xi$.  Let
$\alpha_1 = \Xi_1 (\alpha_0)$.  Then $\alpha_1$ agrees with
$\alpha_0$ on $D_N$, and $u_1 \circ p_* \alpha_1$ is a
polynomial of degree $\leq N$.

\vspace{.1in}

To prove the claim, we construct the power series
$$f = \sum_{i = N + 1 - n_1}^\infty f_i \cdot t^i$$
inductively, coefficient by coefficient.  Let
$(g_1 \circ \gamma)_{n_1}$ be the first non-zero coefficient of
$g_1 \circ \gamma$, and $(u_1 \circ p_* \gamma)^{}_{N+1}$
be the $(N+1)$-st coefficient of $u_1 \circ p_* \gamma$.
We set $$f^{}_{N + 1 - n_1} = - \frac
{(u_1 \circ p_* \gamma)^{}_{N+1}}{(g_1 \circ \gamma)^{}_{n_1}}.$$
This ensures that $(u_1 \circ p_* \alpha_1)^{}_{N+1} = 0$.  Once
$f^{}_{N + 1 - n_1}$ is fixed, the next coefficient,
$f^{}_{N + 2 - n_1}$, is determined uniquely by the requirement
$(u_1 \circ p_* \alpha_1)^{}_{N+2} = 0$, and so on.  Note that
the time-$1$ flow map $\Xi_1$ restricts to an isomorphism of
formal neighborhoods: $\cL (X)_{\alpha_0} \cong \cL (X)_{\alpha_1}$.

The next $d-1$ steps of the construction are completely
analogous.  Each time we obtain an arc $\alpha_k$
($k = 2, \dots, n$) which agrees with $\alpha_{k-1}$ on
$D_N$, and satisfies: $u_i \circ p_* \alpha_k = 
u_i \circ p_* \alpha_{k-1}$, for $i \neq k$, and
$u_k \circ p_* \alpha_k$ is a polynomial of degree $\leq N$.
In the end, we obtain the required arc $\alpha = \alpha_d$.
Claims (i)-(iii) of the lemma follow from the construction,
and the complex analyticity of $\alpha$ follows from Lemma
\ref{prelim3}.
\hfill$\Box$

\vspace{.2in}

\noindent
{\bf Step 3:  General position assumptions.}
Using Lemma \ref{truncate}, we can assume from now on that the
arc $\gamma$ is analytic, and that the image of $\gamma$ is a
local branch of some algebraic curve.  For the arguments in
Steps 5 and 6, we will also need to assume that the complement
$W$ is in general position with respect to $\gamma$.  The purpose
of Step 3 is to specify this general position assumption.

Given any $x_1 \neq x_2$ in $X$, denote by $L(x_1, x_2) \subset U$
the straight line through $x_1$ and $x_2$.  For $x \in X$, let
$$K(x) = \bigcup_{x' \in X \setminus \{ x \}} L(x, x') \subset U,$$
and let $\bar K(x)$ be the closure of $K(x)$ in $U$.  Given an
analytic arc $\alpha \in \cL (X)$, let ${\displaystyle K (\alpha) =
\lim_{t \to 0} \bar K (\alpha(t))} \subset U$.  Note that
$K (\alpha)$ automatically contains $K (\alpha(0))$.

\begin{lemma}\label{cone}
The set $K (\gamma)$ is a closed algebraic cone, with
$\dim K (\alpha) \leq d+1$.
\hfill$\Box$
\end{lemma}

Given an analytic arc $\alpha \in \cL (X)$, we denote by
$\cD_\alpha \subset \C$ a small disc around the origin, such
that $\alpha$ can be viewed as a map $\cD_\alpha \to X$.  We also
let $\cD^\circ_\alpha = \cD_\alpha \setminus \{ 0 \}$.

\begin{lemma}\label{generic}
For a suitably generic choice of the complement $W$ to $V$ in $U$,
there is a positive integer $N = N(W)$, such that for any analytic
$\alpha \in \cA (X, \gamma, N)$, there is a neighborhood $\cU$ of
$0$ in $U$, such that:

(i)   $\dim (W \cap K(\alpha)) \leq 1$.

(ii)  $W \cap X \cap \cU = \{ 0 \}$.

(iii) The intersection $p^{-1} (p \circ \alpha (\cD^\circ_\alpha))
\cap X \cap \cU$ in contained in $X^\circ$.

(iv)  The intersection in part (iii) is transverse.  That is,
$p^{-1} (p \circ \alpha (\cD^\circ_\alpha)) \cap \cU$ is transverse
to $X^\circ \cap U$ as (non locally closed) submanifolds of $\cU$.
\end{lemma}

\noindent
{\bf Proof:}  This is a standard general position argument.
The complement $W$ is chosen as follows.  Let $C(X) \subset U$ be
the tangent cone of $X$ at zero; it is a conical subvariety of $U$
of dimension $d$.  Write $X^{sing} = X \setminus X^\circ$.
For $x \in X$, let $$S(x) = \bigcup_{x' \in X^{sing} \setminus
\{ x \}} L(x, x') \subset U,$$ and let $\bar S(x)$ be the closure
of $S(x)$ in $U$.  Let $S (\gamma) = \lim_{t \to 0} \bar S
(\alpha(t))$.  By analogy with Lemma \ref{cone}, $S(\gamma)$ is a
closed algebraic cone, with $\dim S(\gamma) \leq d$.

Let $c = \dim U - \dim V$.  Define $\cP^\circ$ to be the set of all
affine $c$-planes in $U$ which have a point of non-transverse
intersection with $X^\circ$.  Let $\cP$ be the closure of
$\cP^\circ$ in the Grassmannian of all affine $c$-planes.  For any
$x \in X$, let $\cP(x)$ the set of all $P \in \cP$ passing through
$x$.  By standard general position, $\cP (x)$ is a proper, closed
subvariety of the Grassmannian of all $c$-planes through $x$.  Set
$\cP (\gamma) = \lim_{t \to 0} \cP (\gamma (t))$; it is a proper,
closed subvariety of the Grassmannian of linear $c$-planes in $U$.
To satisfy conditions (i) - (iv) of the lemma, it is enough to
choose $W$ so that:

(1) $\dim (W \cap K(\gamma)) \leq 1$;

(2) $W \cap C(X) = \{ 0 \}$;

(3) $W \cap S(\gamma) = \{ 0 \}$;

(4) $W \notin \cP (\gamma)$.

\noindent
Manifestly, each of these four conditions is Zariski open and
dense.
\hfill$\Box$

\vspace{.1in}

Using Lemmas \ref{generic} and \ref{truncate}, we assume from now
on that the arc $\gamma$ and the complement $W$ satisfy conditions
(i) - (iv) of Lemma \ref{generic}, and that $p\gamma \in
\cL(V)^{\leq N}$ for some large $N \in \Z_+$.

\vspace{.2in}

\noindent
{\bf Step 4: A product decomposition.}

\begin{lemma}\label{cl_im}
The map $p_* : \cL(X) \to \cL(V)$ induces a closed immersion
$$\cL(X)_\gamma \hookrightarrow \cL(V)_{p\gamma}.$$
The differential $d_\gamma p_* : T_\gamma \cL(X) \to
T_{p\gamma} \cL(V)$ is an isomorphism.
\end{lemma}

\noindent
{\bf Proof:}
This is an exercise in the inverse function theorem.
\hfill$\Box$

\vspace{.1in}

Decomposition (\ref{decomp}) induces a decomposition of
formal neighborhoods
$$\cL(V)_{p\gamma} = \cL(V)^{\leq N}_{p\gamma} \times
\cL(V)^{>N}_{p\gamma}.$$
Let $p_* \cL(X)_\gamma \subset \cL(V)_{p\gamma}$ be the
image of the closed immersion in Lemma \ref{cl_im}.
Set $F = p_* \cL(X)_\gamma \cap \cL(V)^{\leq N}_{p\gamma}$.

\begin{lemma}\label{product}
We have $p_* \cL(X)_\gamma \supset \cL(V)^{>N}_{p\gamma}$, and
$$p_* \cL(X)_\gamma \cong F \times \cL(V)^{>N}_{p\gamma}.$$
\end{lemma}

\noindent
{\bf Proof:}
Containment $p_* \cL(X)_\gamma \supset
\cL(V)^{>N}_{p\gamma}$ follows from Lemma \ref{prelim3}.
The proof of the product decomposition is similar to
the proof of Lemma \ref{truncate}.  We sketch it below,
continuing with the notation of the proof of Lemma
\ref{truncate}.

Define an index set $J \subset \{1, \dots, d \} \times \Z_+$
by $J = \{ (i, k) \, | \, k > N - n_i \}$.  Write $\C^J$ for
the set of all maps $a : J \to \C$.  For any $a \in \C^J$, let
$$\xi_t^a = \sum_{(i,k) \in J} a(i,k) \, t^k \cdot \eta_i \, ;$$
this gives a family of vector fields on $X$ parametrized by $D$.
Let $\hat \xi^a$ be the corresponding vector field on $\cL (X)$,
and $\Xi_1^a : \cL (X)^{\xi^a} \to \cL (X)$ be the time-$1$ flow
of $\hat \xi^a$.  Define a map $\Psi : F \times \C^J \to \cL(V)$
by $\Psi (\alpha, a) = p_* \Xi_1^a (\alpha)$, where $\alpha \in
F(A)$ for some test-ring $A$.  It is not hard to check that
$\Psi$ induces an isomorphism $\Psi_{(\gamma, 0)} : F \times D^J
\cong p_* \cL(X)_\gamma$, and that the image $\Psi_{(\gamma, 0)}
(\{ \gamma \} \times D^J) = \cL(V)^{>N}_{p\gamma}$.
\hfill $\Box$

\vspace{.2in}

\noindent
{\bf Step 5: Generic projections.}
In order to prove Theorem \ref{main}, we now need to
identify the finite-dimensional piece $F \subset
\cL(V)^{\leq N}_{p\gamma}$.  Our strategy is to present
$F$ as an intersection of finitely many closed subschemes
of $\cL(V)^{\leq N}_{p\gamma}$; then to analyze each of them
using the results of Section 2.

For $l \in W^*$, define $\pi_l : U \to V \oplus \C$ by
$(v, w) \mapsto (v, l(w))$.  Let $W^{*, \circ} \subset W^*$ be
the set of all $l$, such that $\gamma(t)$ lands in the domain of
injectivity of the restriction ${\pi_l}|^{}_X$ (that is,
$X \cap \pi_l^{-1} (\pi_l \circ \gamma (t)) = \{ \gamma (t) \}$),
for small $t$.

\begin{lemma}\label{gen_proj}
The subset $W^{*, \circ} \subset W^*$ is non-empty and Zariski
open.
\end{lemma}

\noindent
{\bf Proof:}
This follow from condition (i) of Lemma \ref{generic}.
\hfill$\Box$

\vspace{.1in}

Pick a basis $L$ of $W^*$, such that $L \subset W^{*, \circ}$.
For each $l \in L$, let $X_l = \overline{\pi_l (X)} \subset
V \oplus \C$, and $\gamma_l  = \pi_l \circ \gamma :
D \to X_l$.  We can make all the constructions of
Step 4 for each of the arcs $\gamma_l$.  In particular,
we obtain a finite piece $F_l = p_* \cL(X_l)_{\gamma_l}
\cap \cL(V)^{\leq N}_{p\gamma}$, for each $l$.

\begin{lemma}\label{intersect}
We have $F = \bigcap_{l \in L} F_l,$
as subschemes of $\cL(V)^{\leq N}_{p\gamma}$.
\end{lemma}

\noindent
{\bf Proof:}
Let $A$ be a test-ring with maximal ideal $\fm$, and
$\alpha \in \cL(V)^{\leq N}_{p\gamma} (A)$ be a map
${\rm Spec} \, A \times D \to V$.  Let
$D^\circ = {\rm Spec} \, \C ((t)) \subset D$ be the
punctured formal disc, $\gamma^\circ$ be the restriction of
$\gamma$ to $D^\circ$, and $\alpha^\circ$ be the restriction
of $\alpha$ to ${\rm Spec} \, A \times D^\circ$.  Then
there is a unique map $\tilde\alpha^\circ : {\rm Spec} \, A
\times D^\circ \to X$ such that $\tilde\alpha^\circ$
restricted to $D^\circ (= {\rm Spec} \, A/\fm \times D^\circ)$
equals $\gamma^\circ$, and $p \circ \tilde\alpha^\circ =
\alpha^\circ$.

By definition, $\alpha \in F(A)$ if and only if
$\tilde\alpha^\circ$ extends to a map of
${\rm Spec} \, A \times D$.  This will happen if and only if
$l \circ \tilde\alpha^\circ \in A[[t]],$ for all $l \in L$.
But for each individual $l$, this condition is equivalent to
saying that $\alpha \in F_l(A)$.
\hfill$\Box$

\vspace{.2in}

\noindent
{\bf Step 6: Reduction to plane curves.}
We now proceed to identify each of the $F_l$ ($l \in L$).
Consider coordinates $(u_1, \dots, u_d, z)$ on $V \oplus \C$,
such that $\C$ is the $z$-axis, and $V$ is the plane $\{ z =
0 \}$.  Let $\phi_l (u_1, \dots, u_d, z)$ be the defining
equation of the hypersurface $X_l \subset V \oplus \C$.  Define
an evaluation map $\psi_l : \cL(V)^{\leq N} \to \C[x,y]$ by
$$\psi_l: \beta \mapsto \phi_l ((u_1 \circ \beta) (x), \dots,
(u_d \circ \beta) (x), y),$$
and let $f_l = \psi_l (p\gamma)$.

\begin{lemma}
The polynomial $f_l \in \C [x,y]$ has no critical points in
some punctured neighborhood of the origin.
\end{lemma}

\noindent
{\bf Proof:}  This follows from conditions (ii) - (iv) of Lemma
\ref{generic}.
\hfill$\Box$

\vspace{.1in}

The rest of the proof is based on applying Lemma
\ref{plane_curves} to $f = f_l$.  Let $C \subset \C^2$ be the
curve $\{ f (x, y) = 0 \}$, and $(C_1, 0) \subset (C, 0)$ be the
germ of the graph $\{ y = z \circ \gamma_l (x) \}$.  Fix an
integer $M \gg N \cdot \deg (\phi_l)$, and let $P = P_l$ be a
small analytic neighborhood of $f$ in $\C[x,y]^{\leq M}$.  Based
on this data, we can define the locus $R = R_l \subset P_l$
as in Section 2.  Consider the formal neighborhoods
$(R_l)^{}_f \subset (P_l)^{}_f$, and the map
$(\psi_l)^{}_{p\gamma} : \cL(V)^{\leq N}_{p\gamma} \to
(P_l)^{}_f$ induced by $\psi_l$.

\begin{lemma}\label{last_step}
We have: $F_l = (\psi_l)_{p\gamma}^{-1} \, (R_l)^{}_f$.
\end{lemma}

\noindent
{\bf Proof:}  Begin with the inclusion $F_l \subset
(\psi_l)_f^{-1} (R_l)^{}_f$.  Let $A$ be a test-ring, and
$\alpha \in F_l(A)$.  This means that $\alpha \in
\cL(V)^{\leq N}_{p\gamma} (A)$ admits a lift
$\tilde\alpha : {\rm Spec} \, A \times D \to X_l$
such that $\tilde\alpha|^{}_D = \gamma_l$, and
$p \circ \tilde\alpha = \alpha$.  Consider the image
$$\tilde f = (\psi_l)^{}_{p\gamma} \, (\alpha) \in
(P_l)^{}_f \, (A) \subset A[x,y]^{\leq M}.$$
Let $\tilde h \in A [[t]]$ be the composition
$z \circ \tilde\alpha : {\rm Spec} \, A \times D \to \C$.
It is a tautology that $\tilde f$ and $\tilde h$ satisfy
the conditions of Lemma \ref{plane_curves} (iv).  Inclusion
$\subset$ follows.  The opposite inclusion in analogous.
\hfill$\Box$

\vspace{.1in}

To complete the proof of Theorem \ref{main}, we need a
basic result of M. Artin.

\begin{lemma}\label{standard} {\em [Ar, Corollary 2.1]}
Let $E \cong \C^n$ be an affine space, $Z \subset E$ be a closed
algebraic subvariety, and $e$ be a point in $Z$.  Suppose the
analytic germ $(Z, e)$ is reducible: $(Z, e) = (Z_1, e) \cup
(Z_2, e)$.  Then there exists an etale neighborhood
$\pi : \tilde E^\circ \to E^\circ \subset E$ of $e$ in $E$,
and a point $\tilde e \in \pi^{-1} (e)$, such that the preimage
$\tilde Z = \pi^{-1} (Z)$ is reducible as an algebraic
variety: $\tilde Z = \tilde Z_1 \cup \tilde Z_2$, and $\pi$
induces an isomorphism of analytic germs $(\tilde Z_1, \tilde e)
\cong (Z_1, e)$. \hfill $\Box$
\end{lemma}

Let now $E = \prod_{l \in L} \C[x,y]^{\leq M}$, where
$M \gg N \cdot \deg (\phi_l)$ for all $l$.  We have a map
$\psi: \cL(V)^{\leq N} \to E$ given, in components, by the
$\{\psi_l\}_{l \in L}$.  Let $e = \psi(\gamma) =
\{f_l\}_{l \in L}$.  For each $l \in L$, let $Z_l \subset
\C[x,y]^{\leq M}$ be a closed subvariety having $R_l$ as an
analytic irreducible component near $f_l$ (cf. Lemma
\ref{plane_curves} (ii)).  Let $Z = \prod_{l \in L} Z_l \subset E$,
and $Z_1 = \prod_{l \in L} R_l$.  Then $Z_1$ is
an analytic component of $Z$ near $e$.  Applying Lemma
\ref{standard}, we obtain an etale neighborhood
$\pi : \tilde E^\circ \to E^\circ \subset E$ of $e$, a point
$\tilde e \in \pi^{-1} (e)$, and an irreducible component
$\tilde Z_1$ of $\tilde Z = \pi^{-1} (Z)$, such that
$(\tilde Z_1, \tilde e) \cong (Z_1, e)$.  The scheme
$Y$ of Theorem \ref{main} is defined as the fiber product
${\displaystyle Y = \cL(V)^{\leq N} \times^{}_E \tilde Z_1}$, and
the point $y$ is taken to be $(\gamma, \tilde e)$.  The theorem
follows from Lemmas \ref{product}, \ref{intersect}, and
\ref{last_step}.

\section{Examples}

The proof of Theorem \ref{main} described above is far from
being constructive.  Here are a few simple examples where we
can explicitly present the versal deformation of an arc.

\vspace{.1in}

\noindent
{\bf Example 1.}  Let $X = \{ (x,y,z) \in \C^3\; | \; xy = z^2 \}$,
and $\gamma \in \cL^\circ (X)$ be the arc $\gamma (t) = (t, 0, 0)$.
Then the scheme $Y = Y(\gamma)$ of Theorem \ref{main} is just a
double point: $Y = {\rm Spec} \, \C[a]/(a^2)$.  The point
$y \in Y(\C)$ is the only $\C$-point of $Y$, and a versal
deformation $\alpha : Y_y \times D \to X$ of $\gamma$ is given
by $\alpha(a,t) = (t,0,a)$.  Note that this means that 
``topologically'' $\gamma$ is a smooth point of $\cL^\circ(X)$.

\vspace{.1in}

\noindent
{\bf Example 2.}  Let now $X = \{ (x,y,z) \in \C^{2+r}\; | \;
xy = z^2 \}$ be the quadric cone in $\C^{2+r}$, for any $r \geq 1$.
Here $z \in \C^r$, and $z^2 \in \C$ is the standard dot-product
$z \cdot z$.  Let $\gamma \in \cL^\circ (X)$ be the arc
$\gamma (t) = (t, 0, 0)$.  Then $Y = Y(\gamma)$ is the quadric cone
in $\C^r$: $Y = \{ w \in \C^r \; | \; w^2 = 0 \}$ (in its reduced
structure, if $r > 1$).  The point $y \in Y(\C)$ is zero, and a
versal deformation $\alpha : Y_y \times D \to X$ of $\gamma$ is
given by $\alpha(w,t) = (t,0,w)$.  Note that this means that $\gamma$
``lies on a stratum of codimension $r-1$ in $\cL^\circ (X)$.''

\vspace{.1in}

\noindent
{\bf Example 3.}  Let $X$ be as in Example 2, and $\gamma \in
\cL^\circ (X)$ be the arc $\gamma (t) = (t^2, 0, 0)$.  Then
$Y = \{ (a,b,v,w) \in \C \times \C \times\C^r \times \C^r\; | \;
a \, w^2 = v^2, \; b \, w^2 = 2vw \}$ (in its reduced structure).
The point $y \in Y(\C)$ is zero, and a versal deformation
$\alpha : Y_y \times D \to X$ of $\gamma$ is given by
$\alpha(a,b,v,w,t) = (a + b \, t + t^2, w^2, v + t \, w)$.  This
means that $\gamma$ ``lies on a stratum of codimension $2r$ in
$\cL^\circ (X)$.''

\vspace{.1in}

\noindent
MIT, Department of Mathematics, 77 Massachusetts Ave., Room 2-181,
Cambridge, MA 02139

\noindent
{\it grinberg@math.mit.edu}

\vspace{.1in}

\noindent
Harvard University, Department of Mathematics, One Oxford Street,
Cambridge, MA 02138

\noindent
{\it kazhdan@math.harvard.edu}

\end{document}